\documentclass{article}
\usepackage{mathrsfs}[12.5pt]
\usepackage{bm, amscd,  amsfonts, amssymb, cite,texdraw}
\usepackage{amsmath, epsfig, cite, color, colordvi}
\usepackage{latexsym,graphicx}
\usepackage{caption}

\usepackage{caption}

\makeatletter \@addtoreset{equation}{section} \makeatother

\newtheorem{theo}{Theorem}[section]

\newtheorem{coro}[theo]{Corollary}

\newtheorem{lem}[theo]{Lemma}

\def\0{\mathbb 0}

\def\b{\boldsymbol{b}}

\def\x{\mathbf x}
\def\y{\mathbf y}
\def\z{\mathbf z}
\def\u{\mathbf u}
\def\v{\mathbf {v}}

\def\des{\mbox{des}}
\def\asc{\mbox{asc}}
\def\plat{\mbox{plat}}

\def\mid{{\,|\,}}

\def\qed{\hfill \rule{4pt}{7pt}}
\begin{document}

\begin{center}
{\Large  Context-free Grammars and Multivariate
 Stable
 \vskip 2mm

 Polynomials over Stirling Permutations         }

\vskip 3mm

 William Y.C. Chen$^1$,
 Robert X.J. Hao$^2$ and Harold R.L. Yang$^3$\\[5pt]
Center for Combinatorics, LPMC-TJKLC \\
Nankai University, Tianjin 300071, P. R. China

\vskip 3mm

 E-mail: $^1$chen@nankai.edu.cn,
$^2$nalanxindao@163.com, $^3$yangruilong@mail.nankai.edu.cn

\end{center}

\begin{abstract}
Recently, Haglund and Visontai  established the stability of the multivariate Eulerian polynomials
as the generating polynomials of the Stirling permutations, which serves as a
 unification of some results of B\'{o}na, Brenti, Janson, Kuba, and Panholzer concerning Stirling permutations.
Let $B_n(x)$ be the generating polynomials of the descent statistic over Legendre-Stirling permutations, and let
$T_n(x)=2^nC_n(x/2)$, where $C_n(x)$ are the second-order Eulerian polynomials.
 Haglund and Visontai  proposed the problems of finding multivariate stable refinements of the polynomials $B_n(x)$ and $T_n(x)$.
 We obtain
 context-free grammars leading to multivariate stable refinements
 of the polynomials $B_n(x)$ and $T_n(x)$.
 Moreover, the grammars enable us to
 obtain combinatorial interpretations of the multivariate polynomials
 in terms of Legendre-Stirling permutations and marked
  Stirling permutations. Such stable  multivariate
  polynomials provide solutions to two problems posed by Haglund and Visontai.

\end{abstract}



\begin{section}{Introduction}

This paper presents an approach to the construction of stable
combinatorial polynomials from the perspective
of context-free grammars. The framework of using
context-free grammars to generate combinatorial polynomials
was proposed by Chen \cite{Chen1993}.  More specifically,
we introduce the structure of marked Stirling permutations, and
we find context-free grammars that lead to
multivariate stable polynomials over marked Stirling
permutations and Legendre-Stirling permutations. These multivariate
stable polynomials provide solutions to two problems posed by
Haglund and Visontai \cite{HV} in their study of multivariate stable refinements of the
second-order Eulerian polynomials.

Let us first review some backgrounds on the second-order Eulerian polynomials. These polynomials were first introduced by Gessel and Stanley \cite{GS}, which are defined as the generating functions of the descent statistic over Stirling permutations. Recall that a Stirling permutation of order $n$ is a permutation $\pi=\pi_1\pi_2\cdots \pi_{2n-1}\pi_{2n}$ of the multiset $\{1,1,2,2,\ldots,n,n\}$, denoted by $[n]_2$, which satisfies the following condition:  if $\pi_i=\pi_j$ then $\pi_k>\pi_i$ whenever $i<k<j$. For $1\leq i\leq 2n$, we say that $i$ is a descent of $\pi$ if $i=2n$ or $\pi_{i}>\pi_{i+1}$. Analogously, $i$ is called an ascent of $\pi$ if $i=1$ or $\pi_{i-1}<\pi_{i}$. Let  ${Q}_n$ denote the set of Stirling permutations of order $n$.
Let $C(n,k)$ be the number of Stirling permutations of $[n]_2$ with $k$ descents,
and let
$$
C_n(x)=\sum_{k=1}^nC(n,k)x^k.
$$
Gessel and Stanley \cite{GS} showed that
$$
 \sum_{n=0}^\infty S(n+k,k) x^n={\frac{C_n(x)}{(1-x)^{2k+1}}},
 $$
where $S(n,k)$,  as usual, denotes the Stirling number of the second kind.
The numbers $C(n,k)$ are called the second-order Eulerian numbers by Graham, Knuth and Patashnik \cite{GKP}, and accordingly the polynomials $C_n(x)$ are called the second-order Eulerian polynomials by Haglund and Visontai \cite{HV}.

The Stirling permutations were further studied by B\'{o}na \cite{zero}, Brenti \cite{Brenti}, Janson\cite{Janson} and Janson, Kuba and
Panholzer \cite{JKP2011}. B\'{o}na \cite{zero} introduced a statistic, called plateau, on stirling permutations, and proved that ascents, descents and plateaux have the same distribution over $Q_n$.
Given a Stirling permutation $\pi=\pi_1\pi_2\ldots\pi_{2n}\in Q_n$, the index $i$ is called a plateau of $\pi$ if $\pi_{i-1}=\pi_{i}$.
Analogous to that of the classical Eulerian polynomials, B\'{o}na \cite{zero} obtained the real-rootedness of the second-order Eulerian polynomials $C_n(x)$.

\begin{theo}\label{Bona}
For any positive integer $n$, the roots of the polynomial $C_n(x)$ are all real, distinct, and non-positive.
\end{theo}

It should be noted that the real-rootedness of $C_n(x)$ is essentially the real rootedness of the generating function of generalized Stirling permutations obtained by  Brenti \cite{Brenti}.
A permutation $\pi$ of the multiset $\{1^{r_1},2^{r_2},\ldots,n^{r_n}\}$ is called a generalized Stirling permutation of rank $n$ if $\pi$
satisfies the same condition as for a Stirling permutation. Let $Q_n^*$ denote the set of generalized Stirling permutations of rank $n$. In particular, if $r_1=r_2=\cdots=r_n=r$ for some $r$, then $\pi$ is called an $r$-Stirling permutation of order $n$. Let $Q_n(r)$ denote the set of $r$-Stirling permutations of order $n$. It is clear that $1$-Stirling permutations are  ordinary permutations and
$2$-Stirling permutations are the Stirling permutations. Brenti \cite{Brenti} showed that the descent generating polynomials over  $Q_n^*$  have only real roots.

 Janson \cite{Janson} defined  the following trivariate generating function
\begin{align*}
C_n(x,y,z)=\sum_{\pi\in Q_n} x^{\des(\pi)}y^{\asc(\pi)}z^{\plat(\pi)},
\end{align*}
where $\des(\pi)$, $\asc(\pi)$, and $\plat(\pi)$ denote the numbers of descents, the number of ascents, and the number of plateaux of $\pi$, respectively, and proved that $C_n(x,y,z)$ is symmetric in $x,y,z$. This implies the equidistribution  of these three statistics derived by B\'{o}na.

The symmetric property of $C_n(x,y,z)$ was further extended to $r$-Stirling permutations by Janson, Kuba and Panholzer \cite{JKP2011}. For an $r$-Stirling permutation, they introduced the  notion of a $j$-plateau.   For an $r$-Stirling permutation $\pi=\pi_1\pi_2\ldots\pi_{nr}$ and an integer $1\leqslant j\leqslant r-1$, a number $1\leqslant i< nr$ is called a $j$-plateau of $\pi$
if
$\pi_{i}=\pi_{i+1}$ and there are $j-1$ indices $l<i$ such that $\pi_{l}=\pi_i$, i.e., the number $\pi_i$ appears $j$ times
 up to the $i$-th position of $\pi$. Let $j$-plat($\pi$) denote the number of $j$-plateaux of $\pi$. Meanwhile, define a descent and an ascent of $\pi$ similar as ordinary permutations, and let $\des(\pi)$ and $\asc(\pi)$ denote the number of descents and ascents of $\pi$.
Janson, Kuba and Panholzer  \cite{JKP2011} showed that the distribution of ($\des$, $1$-$\plat$, $2$-$\plat$, $\ldots$, $(r-1)$-$\plat$, $\asc$) is symmetric over the set of $r$-Stirling permutations.

Based on the theory of multivariate stable polynomials  recently developed by Borcea and Br\"{a}nd\'{e}n  \cite{BB1, BB, BB2},
Haglund and Visontai \cite{HV} presented a unified approach
 to the stability of the generating functions of Stirling permutations and $r$-Stirling permutations.  A polynomial $f(\z)\in \mathbb{C}[\z]=\mathbb{C}[z_1,z_2,\ldots,z_m]$ is said to be stable, if whenever the imaginary part $ \mathrm {Im}(z_i)>0$ for all $i$ then $f(\mathbf{z})\neq 0$. Clearly, a univariate polynomial $f(z)\in \mathbb{R}[z]$ has only real roots if and only if it is stable.

For the case of univariate real polynomials, P\'{o}lya and Schur \cite{PS} characterized all diagonal operators preserving stability or real-rootedness. Recently, Borcea and Br\"{a}nd\'{e}n \cite{BB1, BB, BB2}  characterized all linear operators preserving stability of multivariate polynomials, see also the survey of Wagner \cite{Wagner}. This implies a characterization
 of linear operators preserving stability of univariate polynomials.

A multivariate polynomial is called multiaffine if the degree
 of each variable is at most $1$. Borcea and Br\"{a}nd\'{e}n showed that each of the operators preserving stability of multiaffine polynomials has a simple form. Using this property, Haglund and Visontai \cite{HV} obtained
  a stable multiaffine refinement of the second-order Eulerian polynomial $C_n(x)$. Given a Stirling permutation $\pi=\pi_1\pi_2\cdots\pi_{2n}\in Q_n$, let
\begin{eqnarray*}
{A}(\pi)&=&\{i|\pi_{i-1}<\pi_{i}\},\\[5pt]
{D}(\pi)&=&\{i|\pi_{i}>\pi_{i+1} \}, \\[5pt]
{P}(\pi)&=&\{i|\pi_{i-1}=\pi_{i} \}
\end{eqnarray*}
denote the set of ascents, the set of descents and  the set plateaux of $\pi$, respectively. Define
$$
C_n(\x,\y,\z)=\sum_{\pi\in{Q}_n}\prod_{i\in {D}(\pi)}x_{\pi_i}\prod_{i\in {A}(\pi)}y_{\pi_{i}}\prod_{i\in {P}(\pi)}z_{\pi_i}.
$$
Haglund and Visontai\cite{HV} proved the stability of $C_n(\x, \y, \z)$.

\begin{theo}\label{cn}
The polynomial
$C_n(\x,\y,\z)$ is stable.
\end{theo}

It is worth mentioning that, as observed by Haglund and Visontai, the recurrence relation between $C_{n-1}(\x,\y,\z)$ and $C_n(\x,\y,\z)$ can   be used to derive the symmetry of $C_n(\x,\y,\z)$, which implies the symmetry of $C_n(x,y,z)$ obtained by Janson, Kuba and Panholzer \cite{JKP2011}.

Moreover, Haglund and Visontai extended the  stability of $C_n(\x,\y,\z)$ to generating polynomials of $r$-Stirling permutations by taking the $j$-plateau statistic into consideration.
Let $ {P}_j(\pi)$ denote the set of $j$-plateaux of $\pi$.
Haglund and Visontai \cite{HV} obtained the following multivariate stable polynomial over $r$-Stirling permutations
\[
E_n(\x,\y,\z_1,\ldots,\z_{r-1})=\sum_{\pi\in {Q}_n(r)}\left(\prod_{i\in {D}(\pi)}x_{\pi_i}\right)\left(\prod_{i\in {A}(\pi)}y_{\pi_i}\right)\prod_{j=1}^{r-1}\left(\prod_{i\in {P}_j(\pi)}z_{j,\pi_i}\right).
\]
They also obtained a similar multivariate stable polynomial for generalized Stirling permutations.

In view of the real-rootedness of $C_n(x)$ and its multivariate stable refinement $C_n(\x,\y,\z)$, Haglund and Visontai posed
 the problem of finding  multivariate stable polynomials as  refinements of the generating polynomials of the descent statistic over Legendre-Stirling permutations. The Legendre-Stirling permutations were introduced
by Egge\cite{Egge} as a generalization of
 Stirling permutations in the study of  Legendre-Stirling numbers of the second kind. For any $n\geq 1$, let $M_n$ be the multiset $\{1,1, \bar{1}, 2,2,\bar{2},\cdots,n,n,\bar{n}\}$. A
permutation $\pi=\pi_1\pi_2\ldots\pi_{3n}$ on $M_n$ is called a Legendre-Stirling   permutation if whenever  $i < j < k$ and $\pi_i=\pi_k$ are both unbarred, then $\pi_j > \pi_i$. For a  Legendre-Stirling permutation $\pi$ on $M_n$, we say that
 $i$ is a descent if either $i=3n$ or
 $\pi_i> \pi_{i+1}$. Let $B_{n,k}$ denote the number of Legendre-Stirling permutations of $M_n$ with $k$ descents. Define
 \[
 B_n(x)=\sum_{k=1}^{2n-1}B_{n,k}x^k.
 \]
Egge proved the real-rootedness of $B_n(x)$.

\begin{theo}\label{Egge}
For $n\geq 1$, $B_n(x)$ has distinct, real, non-positive roots.
\end{theo}

In order to derive a multivariate stable refinement of $B_n(x)$, we introduce an approach of generating stable polynomials by a sequence of grammars. Based on the Stirling grammar given by Chen and Fu \cite{CF}, we find a sequence $G_1, G_2, \ldots$ of context-free  grammars   to generate Legendre-Stirling permutations. We show that the formal derivative with respect to $G_n$ preserves stability  by applying Borcea and Br\"{a}nd\'{e}n's characterization of linear operators  preserving stability.
This leads to a multivariate stable refinement  $B_n(\x, \y, \z,\u,  \v)$ of $B_n(x)$.
On the other hand, according to the grammars, we obtain the following combinatorial interpretation
$$
B_n(\x, \y, \z,\u,  \v)=\sum_{\pi}\prod_{i\in X(\pi)}x_{\pi_i}\prod_{i\in Y(\pi)}y_{\pi_i}\prod_{i\in Z(\pi)}z_{\pi_i}\prod_{i\in U(\pi)}u_{\pi_i}\prod_{i\in V(\pi)}v_{\pi_i}.
$$
The real-rootedness of $B_n(x)$ is a consequence of the
 stability of $B_n(\x, \y, \z,\u,  \v)$ by setting  $v_i=y_i=y$ and $x_i=z_i=u_i=1$ for  $0\leq i\leq n$.

Haglund and Visontai also posed the problem of finding  multivariate stable refinements of the polynomials $T_n(x)$,  which are given by

\begin{equation}\label{tn}
T_n(x)=2^nC_n\left(\frac{x}{2}\right)=\sum_{k}2^{n-k}C(n,k)x^k,
\end{equation}
where $C(n,k)$ and $C_n(x)$, as before, denote the second-order Eulerian numbers and the second-order Eulerian polynomials respectively.  The polynomials $T_n(x)$ were introduced by  Riordan\cite{Riordan}.

In light of the relation (\ref{tn}) between $T_n(x)$ and $C_n(x)$, we introduce the structure of marked Stirling permutations and the following multivariate  polynomials
$$T_n(\x,\y,\z)=\sum_{\pi}\prod_{i\in {D}(\pi) }x_{\pi_{i}}\prod_{i\in {A}(\pi) }y_{\pi_i}\prod_{i\in {P}(\pi) }z_{\pi_i},
$$
where $\pi$ ranges over marked Stirling permutations of $[n]_2$.
We shall show that the polynomials $T_n(\x,\y,\z)$ are stable. The polynomial $T_n(x)$ becomes  the specialization of $T_n(\x,\y,\z)$ by setting $x_i=z_i=1$ and $y_i=x$ for $0\leq i\leq n$. This implies that $T_n(x)$ is real-rooted.

This paper is organized as follows. In Section \ref{grammar}, we give an overview of differential operators associated with
context-free grammars. We find context-free grammars
 to generate the polynomials $C_n(\x,\y,\z)$. In Section \ref{LS}, we obtain context-free grammars that lead to the multivariate
   generating polynomials $B_n(\x, \y, \z,\u,\v)$.
In Section \ref{CSP}, we introduce the structure of marked Stirling permutations, and we give context-free grammars to generate the multivariate polynomials $T_n(\x,\y,\z)$.  In Section \ref{theory}, based on Borcea and Br\"{a}nd\'{e}n's characterization of stability preserving linear operators, we
present an approach to proving the stability of polynomials generated by context-free grammars. In particular,
we prove
the stability of multivariate polynomials $B_n(\x,\y,\z,\u,\v)$ and $T_n(\x,\y,\z)$.
\end{section}

\begin{section}{Context-free grammars}\label{grammar}

In this section, we give an overview of the
idea of using context-free grammars $G$ to generate combinatorial
 polynomials and combinatorial structures as developed by Chen \cite{Chen1993}.
 A context-free grammar $G$ over an
 alphabet $A$ is defined to be a set of production rules.
Roughly speaking,
a production rule means
to substitute a letter in the alphabet $A$ by a
 polynomial in $A$ over a field. Given a context-free grammar,
  one may define a formal derivative $D$ as a linear operator
  on polynomials in $A$, where the action of $D$ on a letter
  is defined by the substitution rule of the grammar and the
  action of $D$ on a product of two polynomials $u$ and $v$
  is defined by the Leibnitz rule, that is,
  \[
  D(uv)=D(u)v+uD(v).
  \]

Many combinatorial polynomials can be generated by context-free grammars. Meanwhile, context-free grammars
 can be used to generate combinatorial structures. More precisely,
 one may use a word on an alphabet to label a combinatorial
 structure such that the context-free grammar serves as
 the procedure to recursively generate the combinatorial
 structures. Such a labeling of a combinatorial structure is
 called a grammatical labeling in \cite{CF}.

For example,   the grammar \[ G=\{a \rightarrow
 ab,\, b \rightarrow b\}\] is used in  \cite{Chen1993} to generate the set of partitions of $[n]$
and the Stirling polynomials,
 $$
 S_n(x)=\sum_{i=0}^nS(n,k)x^k,
 $$
 where $S(n,k)$ denotes the Stirling number of the second kind. For a partition $P$, we label a block of $P$ by letter $b$ and label the partition itself by letter $a$, and we define the weight of
  a partition by the product of its labels.
  So a partition $P$ with $k$ blocks has the  weight $w(P)=ab^k$. For example, the partition $\{\{1,2\},\{3\}\}$ is labeled as follows
\[\stackrel{\{1,2\}}{\mbox{$b$}}
\stackrel{\{3\}}{\mbox{$b$}}
\stackrel{}{\mbox{$a$}}.
\]
In the above notation, we write a partition $P=\{P_1,P_2,\ldots,P_k\}$ of $[n]$ in such a way that the blocks are ordered in the
increasing order of  their minimum elements. Moreover, we put the
letter $a$ at the end of the partition.

Using the above grammatical labeling of a partition,
we deduce that
\begin{equation}\label{partition}
 D^n(a)=\sum_{P}w(P)=\sum_{k=1}^nS(n,k)ab^k.
 \end{equation}
Many properties of the Stirling polynomials follow from the
above expression in terms of the differential operator $D$
with respect to the grammar $G$.

Let us explain how the grammar works for the generation of
partitions. For $n=1$, there is one partition of $[1]$, that is, $\{\{1\}\} $, whose label is $ab$. Assume that we have
generated all the partitions of $[n-1]$ by applying the operator
$D^{n-2}$ to $\{\{1\}\} $ with the initial grammatical labeling.

Let us give an example to demonstrate the action of the
differential operator $D$ with respect to the grammar $G$
to a partition of $[n]$ with the aforementioned grammatical
labeling. Consider the following partition of $\{1, 2, 3, 4, 5, 6\}$
\[\stackrel{\{1,3,6\}}{\mbox{$b$}}\stackrel{\{2,5\}}{\mbox{$b$}}
     \stackrel{\{4\}}{\mbox{$b$}}\stackrel{}{\mbox{$a$}}.\]

If we apply the substitution rule to the letter $a$, then
 we get
$ab$ which we rewrite as $ba$, where $a$ still serves as the label of the new partition, and
$b$ stands for a new block $\{7\}$. In this case, we get a partition
\[\stackrel{\{1,3,6\}}{\mbox{$b$}}\stackrel{\{2,5\}}{\mbox{$b$}} \stackrel{\{4\}}{\mbox{$b$}}
\stackrel{\{7\}}{\mbox{$b$}}\stackrel{}{\mbox{$a$}}.
\]

If we apply the substitution rule to the second letter $b$, then
we get $b$. In this case, we insert the element $7$ in the second block, and we are led to the following partition with
consistent grammatical labeling
\[\stackrel{\{1,3,6\}}{\mbox{$b$}}\stackrel{\{2,5,7\}}{\mbox{$b$}} \stackrel{\{4\}}{\mbox{$b$}}
\stackrel{}{\mbox{$a$}}.
\]

Starting with the empty partition with label $a$, we get
\begin{eqnarray*}
 D(a)&=&\stackrel{\{1\}}{\mbox{$b$}}\stackrel{}{\mbox{$a$}} ,\\[5pt]
 D^2(a)&=& \stackrel{\{1\}}{\mbox{$b$}}\stackrel{\{2\}}{\mbox{$b$}}
 \stackrel{}{\mbox{$a$}}+
 \stackrel{\{1,2\}}{\mbox{$b$}}\stackrel{}{\mbox{$a$}},\\[5pt]
 D^3(a)&=&
     \stackrel{\{1\}}{\mbox{$b$}}\stackrel{\{2\}}{\mbox{$b$}}
     \stackrel{\{3\}}{\mbox{$b$}}\stackrel{}{\mbox{$a$}}
    +\stackrel{\{1\}}{\mbox{$b$}}
    \stackrel{\{2,3\}}{\mbox{$b$}}
    \stackrel{}{\mbox{$a$}}+
    \stackrel{\{1,3\}}{\mbox{$b$}}
    \stackrel{\{2\}}{\mbox{$b$}}
    \stackrel{}{\mbox{$a$}} +
     \stackrel{\{1,2\}}{\mbox{$b$}}
    \stackrel{\{3\}}{\mbox{$b$}}
    \stackrel{}{\mbox{$a$}} +
     \stackrel{\{1,2,3\}}{\mbox{$b$}}\stackrel{}{\mbox{$a$}}.\\[5pt]
\end{eqnarray*}

Without considering the combinatorial structures during
the applications of the differential operator $D$, we may directly
compute $D^n(x)$ to derive the Stirling polynomials $S_n(x)$.

As the second example, we consider the   context-free grammar \[G=\{x\rightarrow xy,\, y\rightarrow xy\} \]
 introduced by Dumont \cite{Dumont} which  is used to
 compute
 the Eulerian polynomials. For a permutation $\pi=\pi_1\pi_2\ldots\pi_{n}$ of $[n]$, let
\begin{eqnarray*}
 {A}(\pi)&=&\{i\, |\, \pi_{i-1}<\pi_i\},\\[5pt]
 {D}(\pi)&=&\{i\, |\, \pi_{i}>\pi_{i+1} \}
\end{eqnarray*}
denote the set of ascents and the set of descents of $\pi$, respectively. Here we set $\pi_0=\pi_{n+1}=0$. In other words,
  for any  permutation $\pi$ of $[n] $, $ 1$ is always an ascent and $n$ is always a descent. An element $\pi_i$
 is called  a descent top of $\pi$  if $i\in D(\pi)$, and $\pi_i$
   is called an ascent top if $i\in A(\pi)$, see
   Haglund and Visontai \cite{HV}.

  The grammatical labeling of a permutation $\pi$ is defined as follows. If $\pi_i$ is an ascent top of $\pi$,
    then we label $\pi_{i-1}$ with  the letter $x$.
    If $\pi_i$ is a descent top,
     then we label $\pi_i$ by the letter $y$.
     For this labeling,  the weight of $\pi$ is given by
\[
w(\pi)=x^{|A(\pi)|}y^{|D(\pi)|}.
\]
 Then for $n\geqslant 1$, we have
\[
D^{n}(x)=\sum_{\pi\in S_n}w(\pi)=\sum_{m=1}^{n}A(n,m)y^{m}x^{n+1-m},
\]
where $A(n,m)$ is the Eulerian number, namely, the number of permutations of $[n]$ with $m$ descents, see Dumont \cite{ Dumont}.

For $n=1$, there is only one permutation of $[1]$, that is $1$, whose label is $xy$. Assume that we have
generated all the permutations of $[n-1]$ by applying the operator $D^{n-2}$ to $1$.

Next we give an example to illustrate the action of $D$ on
a permutation of $[6]$. Take a permutation
\[
\stackrel{}{\mbox{$x$}}\ \stackrel{3}{\mbox{$y$}}\ \stackrel{2}{\mbox{$x$}}\ \stackrel{5}{\mbox{$x$}}\ \stackrel{6}{\mbox{$y$}}\ \stackrel{4}{\mbox{$y$}}\ \stackrel{1}{\mbox{$y$}}.
\]
If we apply the substitution rule $x\rightarrow xy$
to the third letter $x$, we insert $7$ after $5$.
As for the grammatical labeling, we keep all the
labels and assign the element $7$ a new label $y$
as if it comes from the substitution rule
$x\rightarrow xy$. Indeed, it is easily checked that
what we get is a permutation with a consistent grammatical labeling, namely,
\[
\stackrel{}{\mbox{$x$}}\ \stackrel{3}{\mbox{$y$}}\ \stackrel{2}{\mbox{$x$}}\ \stackrel{5}{\mbox{$x$}}\ \stackrel{7}{\mbox{$y$}}\ \stackrel{6}{\mbox{$y$}}\ \stackrel{4}{\mbox{$y$}}\ \stackrel{1}{\mbox{$y$}}.
\]

Similarly, if we apply the substitution rule $y\rightarrow xy$ to the second letter $y$, then we insert $7$ after $6$.
 In this case, we need to change the label of $6$ from $y$ into $x$, and assign $y$ to the new element  $7$. In other words, the label $y$ becomes $xy$ just like the substitution rule. So we
 get the following permutation with a grammatical labeling,
\[
\stackrel{}{\mbox{$x$}}\ \stackrel{3}{\mbox{$y$}}\ \stackrel{2}{\mbox{$x$}}\ \stackrel{5}{\mbox{$x$}}\ \stackrel{6}{\mbox{$x$}}\ \stackrel{7}{\mbox{$y$}}\ \stackrel{4}{\mbox{$y$}}\ \stackrel{1}{\mbox{$y$}}.
\]

Indeed, the above examples indicate
 that permutations of $[n]$
and the Eulerian polynomials $A_n(x)$ can be generated by the
operator $D$ associated with the grammar $G$.

In order to generate combinatorial structures
with more parameters, we may use  a sequence of
 grammars. Let us consider the
the multivariate refinement of Eulerian polynomials $A_n(\x,\y)$
introduced by Haglund and Visontai \cite{HV}, which involve
the sets of   ascent tops and descent tops, not just the numbers of ascents and descents. More precisely,
\[
A_n(\x,\y)=\sum_{\pi\in S_n}\prod_{i \in A(\pi)}x_{\pi_i}\prod_{i \in D(\pi)}y_{\pi_i}.
\]

We shall introduce a sequence of grammars $\{G_n\}$ to generate the multivariate polynomials $A_n(\x,\y)$.

For $n\geq 1$, define
$$
G_n=\{x_i\rightarrow x_ny_n,y_i\rightarrow x_ny_n, 0\leq i<n\},
$$
and denote by $D_n$ the formal differential operator with respect to $G_n$. The multivariate polynomials $A_n(\x,\y)$ can be
generated by the sequence of grammars $G_n$.

\begin{theo}\label{Pgrammar} For $n\geq 1$,
we have
$$D_nD_{n-1}\cdots D_{1}(x_0)=A_n(\mathbf{x},\mathbf{y}).$$
\end{theo}

\begin{pf}
We define the grammatical labeling of a permutation $\pi$ as follows.
For a permutation $\pi$, if $\pi_i$ is an ascent top, we label $\pi_{i-1}$ by the letter $x_{\pi_i}$; if $\pi_i$ is a descent top, we label $\pi_i$ by the letter $y_{\pi_i}$. So the weight of $\pi$ is given by
$$
w(\pi)=\prod_{i \in A(\pi)}x_{\pi_i}\prod_{i \in D(\pi)}y_{\pi_i}.
$$

We proceed to show by induction that $D_nD_{n-1}\cdots D_{1}(x_0)$ equals the sum of the weights of   permutations of $[n]$. For $n=1$, the theorem is valid since the weight of the permutation $1$ is $x_1y_1$.
Assume that the theorem holds for $n-1$, that is,
\[
D_{n-1}\cdots D_{1}(x_0)=\sum_{\pi\in S_{n-1}}w(\pi).
\]
We now use an example to illustrate the action of $D_n$ on a permutation of $[n-1]$.
Let $\pi=325641$. The grammatical labeling is as follows
\[
\stackrel{}{\mbox{$x_3$}}\ \stackrel{3}{\mbox{$y_3$}}\ \stackrel{2}{\mbox{$x_5$}}\ \stackrel{5}{\mbox{$x_6$}}\ \stackrel{6}{\mbox{$y_6$}}\ \stackrel{4}{\mbox{$y_4$}}\ \stackrel{1}{\mbox{$y_1$}}.
\]

If we apply the substitution rule $x_6\rightarrow x_7y_7$ to the letter $x_6$, we define the action as
 the insertion of $7$ immediately after $5$.
 The labels of $5$ and $7$ will be changed to
 $x_7$ and $y_7$ as given by the grammar. It is not hard to see that
   the permutation we obtain has a consistent grammatical labeling,
\[
\stackrel{}{\mbox{$x_3$}}\ \stackrel{3}{\mbox{$y_3$}}\ \stackrel{2}{\mbox{$x_5$}}\ \stackrel{5}{\mbox{$x_7$}}\ \stackrel{7}{\mbox{$y_7$}}\ \stackrel{6}{\mbox{$y_6$}}\ \stackrel{4}{\mbox{$y_4$}}\ \stackrel{1}{\mbox{$y_1$}}.
\]
Similarly, if we apply the substitution rule $y_6\rightarrow x_7y_7$ to the letter $y_6$, we obtain a permutation with a consistent
grammatical labeling
\[
\stackrel{}{\mbox{$x_3$}}\ \stackrel{3}{\mbox{$y_3$}}\ \stackrel{2}{\mbox{$x_5$}}\ \stackrel{5}{\mbox{$x_6$}}\ \stackrel{6}{\mbox{$x_7$}}\ \stackrel{7}{\mbox{$y_7$}}\ \stackrel{4}{\mbox{$y_4$}}\ \stackrel{1}{\mbox{$y_1$}}.
\]
 It
is clear that all permutations of $[n]$ can be obtained   this way.
So we conclude that
$$
D_nD_{n-1}\cdots D_{1}(x_0)=D_n\left(\sum_{\pi\in S_{n-1}}w(\pi)\right)=\sum_{\sigma\in S_n}w(\sigma).
$$
Hence the theorem holds for all positive numbers $n$ by induction.
\qed
\end{pf}

For $n=0$, the empty permutation is labeled
by $x_0$. The values of $A_n(\x,\y)$ for $n=1,2,3$ are given below.
\begin{eqnarray*}
D_1(x_0)&=&\stackrel{}{\mbox{$x_1$}}\  \stackrel{1}{\mbox{$y_1$}},\\[5pt]
D_2D_1(x_0)&=&
\stackrel{}{\mbox{$x_2$}}\ \stackrel{2}{\mbox{$y_2$}}\
\stackrel{1}{\mbox{$y_1$}}+
\stackrel{}{\mbox{$x_1$}}\ \stackrel{1}{\mbox{$x_2$}}\
\stackrel{2}{\mbox{$y_2$}},\\[5pt]
D_3D_2D_1(x_0)&=&
\stackrel{}{\mbox{$x_3$}}\  \stackrel{3}{\mbox{$y_3$}}\ \stackrel{2}{\mbox{$y_2$}}\ \stackrel{1}{\mbox{$y_1$}} +
\stackrel{}{\mbox{$x_2$}}\  \stackrel{2}{\mbox{$x_3$}}\
\stackrel{3}{\mbox{$y_3$}}\ \stackrel{1}{\mbox{$y_1$}}+
\stackrel{}{\mbox{$x_2$}}\  \stackrel{2}{\mbox{$y_2$}}\
\stackrel{1}{\mbox{$x_3$}}\ \stackrel{3}{\mbox{$y_3$}} +
\stackrel{}{\mbox{$x_3$}}\ \stackrel{3}{\mbox{$y_3$}}\
\stackrel{1}{\mbox{$x_2$}}\ \stackrel{2}{\mbox{$y_2$}} \\
& & \quad + \stackrel{}{\mbox{$x_1$}}\ \stackrel{1}{\mbox{$x_3$}}\ \stackrel{3}{\mbox{$y_3$}}\ \stackrel{2}{\mbox{$y_2$}}+
\stackrel{}{\mbox{$x_1$}}\ \stackrel{1}{\mbox{$x_2$}}\
\stackrel{2}{\mbox{$x_3$}}\ \stackrel{3}{\mbox{$y_3$}}.
\end{eqnarray*}

 Let us now consider the grammar to generate Stirling permutations.
 Chen and Fu \cite{CF} showed that the grammar
$$
G=\{x\rightarrow x^2y, y\rightarrow x^2y\}
$$
can be used to generate Stirling permutations.
Let $D$ denote the differential operator associated with
the grammar $G$. It has been shown in \cite{CF} that
$$
D^n(x)=x\sum_{m=1}^nC(n,m)x^{2n-m}y^{m},
$$
where $C(n,m)$ denotes the second-order Eulerian number. Notice
  that
\[
D^n(x)\mid_{x=1}=C_n(y),
\]
where $C_n(y)$ is the second-order Eulerian polynomial.

The grammatical labeling of a Stirling permutation is defined as follows. For a Stirling permutation $\pi$, if $i\in D(\pi)$, we label $\pi_i$ by $y$; if $i\in A(\pi)$ or $i\in P(\pi)$, we label $\pi_{i-1}$ by $x$.
For example, the Stirling permutation $\pi=233211$ has the following grammatical labeling
\[
\stackrel{}{\mbox{$x$}}\ \stackrel{2}{\mbox{$x$}}\
\stackrel{3}{\mbox{$x$}}\ \stackrel{3}{\mbox{$y$}}\ \stackrel{2}{\mbox{$y$}}\ \stackrel{1}{\mbox{$x$}}\ \stackrel{1}{\mbox{$y$}}.
\]

Next we show that one can use a refinement of the grammar $G$
to derive the multivariate polynomials $C_n(\x,\y,\z)$ of
Haglund and Visontai \cite{HV}. Recall that
\[
C_n(\mathbf{x},\mathbf{y},\mathbf{z})=\sum_{\pi\in{Q}_n}\prod_{i\in {A}(\pi)}x_{\pi_i}\prod_{i\in {D}(\pi)}y_{\pi_{i}}\prod_{i\in {P}(\pi)}z_{\pi_i}.
\]

As a refinement of the grammar $G$, we define
\[
G_n=\{x_i\rightarrow x_ny_nz_n, y_i\rightarrow x_ny_nz_n, z_i\rightarrow x_ny_nz_n, 0\leq i<n\}.
\]
and we denote by $D_n$ the differential operator
associated with the grammar $G_n$.

\begin{theo}\label{Sgrammar}
For $n\geqslant 1$,
we have
\[
D_nD_{n-1}\cdots D_{1}(z_{0})=C_n(\x,\y,\z).
\]
\end{theo}

\begin{pf}
First, let us define the grammatical labeling of a Stirling permutation $\pi$.
For a Stirling permutation $\pi$, if $\pi_i$ is an ascent top, we label $\pi_{i-1}$ by the letter $x_{\pi_i}$; if $\pi_i$ is a descent top, we label $\pi_{i}$ by the letter $y_{\pi_i}$; and if $\pi_i$ is a plateau, we label $\pi_{i-1}$ by the letter $z_{\pi_i}$. For this labeling,  the weight of $\pi$ is given by
$$
w(\pi)=\prod_{i\in {A}(\pi)}x_{\pi_i}\prod_{i\in {D}(\pi)}y_{\pi_{i}}\prod_{i\in {P}(\pi)}z_{\pi_i}.
$$
We aim to show that $D_nD_{n-1}\cdots D_{1}(z_0)$ equals the sum of weights of Stirling permutations of $[n]_2$.  Let us use
induction on $n$.
 The theorem is obvious for $n=0$ since the weight of the empty permutation is $z_0$. Assume that the theorem holds for $n-1$, that is,
\[
D_{n-1}\cdots D_{1}(z_0)=\sum_{\pi\in Q_{n-1}}w(\pi).
\]

Let us use an example to demonstrate
the action of $D$ on a Stirling permutation of $[n-1]_2$.  Let $\pi=233211$. The  grammatical labeling of $\pi$ is
as follows
\[
\stackrel{}{\mbox{$x_2$}}\ \stackrel{2}{\mbox{$x_3$}}\
\stackrel{3}{\mbox{$z_3$}}\ \stackrel{3}{\mbox{$y_3$}}\ \stackrel{2}{\mbox{$y_2$}}\ \stackrel{1}{\mbox{$z_1$}}\ \stackrel{1}{\mbox{$y_1$}}.
\]
In general, if we apply a substitution rule of $G_4$
 to any letter in $\pi$, we get $x_4y_4z_4$. Here we insert
 the two elements $44$ after the element whose label
 is replaced by the substitution rule,
 and we use the labels $x_4$, $y_4$ and $z_4$ to
 relabel the three elements that are affected by the
 substitution. For example, if we apply the substitution rule $x_2\rightarrow x_4 y_4 z_4$ to the above Stirling
 permutation, then we get a Stirling permutation
with the following grammatical labeling
\[
\stackrel{}{\mbox{$x_4$}}\ \stackrel{4}{\mbox{$z_4$}}\
\stackrel{4}{\mbox{$y_4$}}\ \stackrel{2}{\mbox{$x_3$}}\
\stackrel{3}{\mbox{$z_3$}}\ \stackrel{3}{\mbox{$y_3$}}\ \stackrel{2}{\mbox{$y_2$}}\ \stackrel{1}{\mbox{$z_1$}}\ \stackrel{1}{\mbox{$y_1$}}.
\]
It is easily seen that the application of any
substitution rule of $G_n$ to any Stirling permutation
of $[n-1]_2$ leads to a Stirling permutation of $[n]_2$ with a consistent
grammatical labeling.
Hence we deduce that
$$
D_nD_{n-1}\cdots D_{1}(z_{0})=D_n\left(\sum_{\pi\in Q_{n-1}}w(\pi)\right)=\sum_{\sigma\in Q_{n}}w(\sigma).
$$
Thus, the theorem holds for $n$. This completes the proof.
\qed
\end{pf}

For $n=0$, the empty permutation is labeled by $z_0$. The values of the polynomials $C_n(\x,\y,\z)$ for $n=1,2$ are   as follows,
\begin{eqnarray*}
D_1(z_{0})&=&\stackrel{}{\mbox{$x_1$}}\ \stackrel{1}{\mbox{$z_1$}}\
\stackrel{1}{\mbox{$y_1$}},\\[5pt]
D_2D_1(z_{0})&=&
\stackrel{}{\mbox{$x_2$}}\ \stackrel{2}{\mbox{$z_2$}}\
\stackrel{2}{\mbox{$y_2$}}\ \stackrel{1}{\mbox{$z_1$}}\
\stackrel{1}{\mbox{$y_1$}}+
\stackrel{}{\mbox{$x_1$}}\ \stackrel{1}{\mbox{$x_2$}}\
\stackrel{2}{\mbox{$z_2$}}\ \stackrel{2}{\mbox{$y_2$}}\
\stackrel{1}{\mbox{$y_1$}}+
\stackrel{}{\mbox{$x_1$}}\ \stackrel{1}{\mbox{$z_1$}}\
\stackrel{1}{\mbox{$x_2$}}\ \stackrel{2}{\mbox{$z_2$}}\
\stackrel{2}{\mbox{$y_2$}}.
\end{eqnarray*}

We shall give further refinements of the above two sequences of grammars as solutions to the problems  of Haglund and Visontai \cite{HV}. On one hand, we use
these refined grammars to construct multivariate polynomials for Legendre-Stirling permutations and marked Stirling permutations.
On the other hand, we use the grammars to
construct stability preserving operators
leading to the stability of the multivariate
polynomials.

\end{section}

\begin{section}{ Legendre-Stirling permutations}\label{LS}

In this section, we introduce several statistics
on Legendre-Stirling permutations of
\[M_n=\{1,1, \bar{1}, 2,2,\bar{2},\cdots,n,n,\bar{n}\}.\]
In terms of these statistics, we obtain multivariate
polynomials $B_n(\x,\y,\z,\u,\v)$ as refinements of $B_n(x)$.
 In fact, the combinatorial construction of the
  multivariate polynomials $B_n(\x,\y,\z,\u,\v)$ is
   obtained from  further refinements of the grammars
   to generate permutations and Stirling permutations with respect to the numbers of descents.
   Using these grammars, we derive the combinatorial
   interpretation by giving a suitable grammatical labeling.
 In Section \ref{theory}, we shall use grammars to
prove the stability of $B_n(\x,\y,\z,\u,\v)$.  This leads to
a solution to the problem of Haglund and Visontai.

 Let  $L_n$ denote the set of Legendre-Stirling permutations of $M_n$. For a Legendre-Stirling permutation $\pi\in L_n$, define
\begin{eqnarray*}
 {X}(\pi)&=&\{i\, |\, \pi_{i-1} \leqslant \pi_i, \mbox{$\pi_i$ is unbarred and appears the first time}\},\\[5pt]
 {Y}(\pi)&=&\{i\, |\, \pi_{i} > \pi_{i+1} \mbox{ and $\pi_i$ is unbarred}\},\\[5pt]
 {Z}(\pi)&=&\{i\, |\, \pi_{i-1} \leqslant \pi_i,\mbox{$\pi_i$ is unbarred and appears the second time}\},\\[5pt]
 {U}(\pi)&=&\{i\, |\, \pi_{i-1} \leqslant \pi_i \mbox{ and $\pi_i$ is barred} \},\\[5pt]
 {V}(\pi)&=&\{i\, |\, \pi_{i} > \pi_{i+1} \mbox{ and $\pi_i$ is barred} \}.
\end{eqnarray*}
Here we  set $\pi_0=\pi_{3n+1}=0$.

For example, let $\pi=\bar{1}1\bar{2}2332\bar{3}1$.
 Then we have $X(\pi)=\{2,4,5\}$, $Y(\pi)=\{6,9\} $, $Z(\pi)=\{6\} $, $U(\pi)=\{1,3,8\} $ and $V(\pi)=\{8\} $.

Define
$$
B_n(\x, \y, \z,\u,  \v)=\sum_{\pi\in L_n}\prod_{i\in X(\pi)}x_{\pi_i}\prod_{i\in Y(\pi)}y_{\pi_i}\prod_{i\in Z(\pi)}z_{\pi_i}\prod_{i\in U(\pi)}u_{\pi_i}\prod_{i\in V(\pi)}v_{\pi_i}.
$$
We define the grammars $\{G_{n}\}$ as follows,
\begin{eqnarray*}
 G_{2n-1}&=&\{x_i,y_i,z_i,u_i,v_i\rightarrow u_n v_n,  0\leqslant i<n\},\\[5pt]
 G_{2n}&=&\{x_i,y_i,z_i,u_i,v_i\rightarrow x_ny_nz_n, 0\leqslant i<n;\\[5pt]
 &&\quad \quad u_{n}\rightarrow x_{n}z_nu_{n}, v_{n}\rightarrow x_{n}y_{n}z_n \}.\\
\end{eqnarray*}
Notice that $G_{2n-1}$ is a refinement of the grammar
\[
G=\{x\rightarrow xy,y\rightarrow xy\}.
\]
and $G_{2n}$ is a refinement of the grammar
\[
G=\{x\rightarrow x^2y,y\rightarrow x^2y\}
\]

The grammatical labeling of a Legendre-Stirling permutation is defined
 as follows. Let $\pi$ be a Legendre-Stirling permutation on $M_n$.
 For $i\in X(\pi)$, $i\in Z(\pi)$ or $i\in U(\pi) $, we label $\pi_{i-1}$ by the letter $x_{\pi_i}$, $z_{\pi_i}$ or $u_{\pi_i} $, respectively; for $i\in Y(\pi)$ or $i\in V(\pi)$, we label $\pi_i$ by the letter $y_{\pi_i}$ or $v_{\pi_i}$, respectively. For example, the above Legendre-Stirling permutation
$\pi=\bar{1}1\bar{2}2332\bar{3}1$ has the following grammatical labeling
\[
\stackrel{}{\mbox{$u_1$}}\ \stackrel{\bar{1}}{\mbox{$x_1$}}\ \stackrel{1}{\mbox{$u_2$}}\ \stackrel{\bar{2}}{\mbox{$x_2$}}\ 
\stackrel{2}{\mbox{$x_3$}}\ \stackrel{3}{\mbox{$z_3$}}\ \stackrel{3}{\mbox{$y_3$}}\ \stackrel{2}{\mbox{$u_3$}}\ \stackrel{\bar{3}}{\mbox{$v_3$}}\ \stackrel{1}{\mbox{$y_1$}}.
\]

 The following theorem shows that the polynomials
 $B_n(\x, \y, \z,\u,  \v)$ can be generated by the grammars
 $G_n$.

\begin{theo}\label{Bgrammar}
For $n\geqslant 1$, we have
\begin{equation}
D_{2n}D_{2n-1}\cdots D_{1}(x_0)=B_n(\x, \y, \z,\u,  \v).
\end{equation}
\end{theo}

\begin{pf}
We use induction on $n$. The case for $n=0$ is obvious since the empty permutation is labeled by $x_0$.
Assume that the theorem holds for $n-1$, that is,
\begin{equation}
D_{2n-2}\cdots D_{1}(x_0)=\sum_{\pi\in L_{n-1}}w(\pi).
\end{equation}

Note that any  Legendre-Stirling permutation of $M_n$
 can be obtained from a Legendre-Stirling permutation  of $M_{n-1}$ through two operations: (1) Insert  a barred element $\bar{n}$; (2) Insert  two elements $nn$. We use an example to show that the operators $D_{2n-1}$ and $D_{2n}$ correspond to these two operations.

Consider the Legendre-Stirling permutation $\pi=\bar{1}1\bar{2}2332\bar{3}1$, whose grammatical labeling is given by
\[
\stackrel{}{\mbox{$u_1$}}\ \stackrel{\bar{1}}{\mbox{$x_1$}}\ \stackrel{1}{\mbox{$u_2$}}\ \stackrel{\bar{2}}{\mbox{$x_2$}}\ 
\stackrel{2}{\mbox{$x_3$}}\ \stackrel{3}{\mbox{$z_3$}}\ \stackrel{3}{\mbox{$y_3$}}\ \stackrel{2}{\mbox{$u_3$}}\ \stackrel{\bar{3}}{\mbox{$v_3$}}\ \stackrel{1}{\mbox{$y_1$}}.
\]

The first operation is just the procedure of generating permutations. In general, if we apply a substitution rule
of $G_7$ to $\pi$, we always get $u_4v_4$.
Here we insert $\bar{4}$ after the element
whose label is replaced by the substitution rule. At the same time, we relabel the two involved elements by the letters $u_4$ and $v_4$. For example, if we apply the substitution rule $z_3\rightarrow u_4v_4$ to $\pi$, then we obtain a Legendre-Stirling permutation with a consistent grammatical labeling
\[
\stackrel{}{\mbox{$u_1$}}\ \stackrel{\bar{1}}{\mbox{$x_1$}}\ \stackrel{1}{\mbox{$u_2$}}\ \stackrel{\bar{2}}{\mbox{$x_2$}}\ 
\stackrel{2}{\mbox{$x_3$}}\ \stackrel{3}{\mbox{$u_4$}}\
\stackrel{\bar{4}}{\mbox{$v_4$}}\  \stackrel{3}{\mbox{$y_3$}}\ \stackrel{2}{\mbox{$u_3$}}\ \stackrel{\bar{3}}{\mbox{$v_3$}}\ \stackrel{1}{\mbox{$y_1$}}.
\]

As for the second operation, consider the above Legendre-Stirling permutation $\sigma=\bar{1}1\bar{2}23\bar{4}32\bar{3}1$.
The two substitution rules $u_4\rightarrow x_4z_4u_4$ and $v_4\rightarrow x_4y_4z_4$ of $G_8$ correspond to
 the operations of inserting two elements $44$ before  $\bar{4}$ or after $\bar{4}$, respectively. So we get two Legendre-Stirling permutations
$\bar{1}1\bar{2}2344\bar{4}32\bar{3}1$ or
$\bar{1}1\bar{2}23\bar{4}4432\bar{3}1$.

Next we consider the rest of substitution rules of $G_8$.
 If we apply any of the remaining substitution rules of $G_8$ to $\sigma$, we always get $x_4y_4z_4$. Here we insert two elements $44$ into $\sigma$
between $\pi_i$ and $\pi_{i+1}$, which are elements less than $4$.
For example, by applying the production rule $u_2
\rightarrow x_4y_4z_4$, we obtain the Legendre-Stirling permutation
\[
\stackrel{}{\mbox{$u_1$}}\ \stackrel{\bar{1}}{\mbox{$x_1$}}\ \stackrel{1}{\mbox{$x_4$}}\ \stackrel{4}{\mbox{$z_4$}}\ 
\stackrel{4}{\mbox{$y_4$}}\  \stackrel{\bar{2}}{\mbox{$x_2$}}\ 
\stackrel{2}{\mbox{$x_3$}}\ \stackrel{3}{\mbox{$u_4$}}\
\stackrel{\bar{4}}{\mbox{$v_4$}}\  \stackrel{3}{\mbox{$y_3$}}\ \stackrel{2}{\mbox{$u_3$}}\ \stackrel{\bar{3}}{\mbox{$v_3$}}\ \stackrel{1}{\mbox{$y_1$}}.
\]

It can be checked that any of applications of the
substitution rules of $G_8$ to $\sigma$ leads to consistent grammatical labelings. Moreover, it can be verified that
 the action of $D_{2n}D_{2n-1}$ on the Legendre-Stirling permutations of $M_{n-1}$ generates all the Legendre-Stirling permutations of $M_n$.
 So we conclude that
$$
D_{2n}D_{2n-1}\cdots D_{1}(x_0)=\sum_{\pi\in L_n}w(\pi).
$$
 Then the theorem follows by induction.
 \qed
\end{pf}

For $n=0$, the empty permutation is labeled by $x_0$, and $B_1(\x,\y,\z,\u,\v)$ is calculated as follows,
\begin{eqnarray*}
D_1(x_0)&=&\stackrel{}{\mbox{$u_1$}}\
\stackrel{\bar{1}}{\mbox{$v_1$}},\\[5pt]
D_2D_1(x_0)&=&\stackrel{}{\mbox{$x_1$}}\ \stackrel{1}{\mbox{$z_1$}}\ \stackrel{1}{\mbox{$u_1$}}\ \stackrel{\bar{1}}{\mbox{$v_1$}}+
\stackrel{}{\mbox{$u_1$}}\
\stackrel{\bar{1}}{\mbox{$x_1$}}\
\stackrel{1}{\mbox{$z_1$}}\ \stackrel{1}{\mbox{$y_1$}}.
\end{eqnarray*}

\end{section}

\begin{section}{Marked Stirling permutations}\label{CSP}

In this section, we introduce the structure of marked Stirling permutations, and we define several statistics in
 order to construct   multivariate polynomials $T_n(\x,\y,\z)$ as refinements of $T_n(x)$.  We also give a
 sequence of grammars to generate $T_n(\x,\y,\z)$ as well as
 marked Stirling permutations with suitable
 grammatical labelings. By using the grammars, the stability of $T_n(\x,\y,\z)$ can be established in Section \ref{theory}.
 This gives a solution to the problem of Haglund and Visontai
  concerning a stable refinement of $T_n(x)$.

A marked Stirling permutation is defined by the following marking rule.
Given a Stirling permutation $\pi=\pi_1\pi_2\cdots\pi_{2n}$,
if $\pi_i$  is an element of $\pi$ such that  $\pi_i$ occurs the second time in $\pi$ and  $\pi_i<\pi_{i+1}$, then we may mark
the element $\pi_i$.  We denote a marked element $i$ by $\bar{i}$.
A marked Stirling permutation is a Stirling permutation
with some elements marked according to the above rule.

For example, there is only one
 marked Stirling permutation of $[1]_2$: $11$, whereas
 there are four marked Stirling permutations of $[2]_2$:
\[2211, 1221, 1122,1\bar{1}22.\]

Let $\bar{Q}_n$ denote the set of marked Stirling permutations of $[n]_2$. We use $ {A(\pi)}$, $ {D(\pi)} $, $ {P(\pi)}$ to denote the set of descents, the set of ascents and the set of plateaux of $\pi$. More precisely, given a marked Stirling permutation $\pi=\pi_1\pi_2\ldots\pi_{2n}\in \bar{Q}_n$, let
\begin{eqnarray*}
{A}(\pi)&=&\{i\, |\, \pi_{i-1}<\pi_{i}\},\\[5pt]
{D}(\pi)&=&\{i\, |\, \pi_{i}>\pi_{i+1} \}, \\[5pt]
{P}(\pi)&=&\{i\, |\, \pi_{i-1}=\pi_{i} \}
\end{eqnarray*}
denote the set of ascents, the set of descents and  the set of plateaux of $\pi$, respectively.
Let $T(n,m)$ be the number of marked Stirling permutations of $[n]_2$ with $m$ descents. It follows from relation (\ref{tn}) that
 \[T_n(x)=\sum_{m=1}^nT(n,m)x^m.
\]

Note that Riordan \cite{Riordan} introduced the polynomials
$T_n(x)$ and proved that $T_n(1)$  equals the Schr\"{o}der number,
namely, the number of series-reduced rooted trees with $n+1$ labeled leaves.

We shall prove that the polynomials $T_n(x)$ can be
generated by the  grammar $G$ defined by
\[G=\{x\rightarrow x^2y, y\rightarrow 2x^2y\}.\]
The grammatical labeling of a marked Stirling permutation
can be described as follows. Let $\pi$ be a marked Stirling permutation of $[n]_2$. If $i\in D(\pi)$, we label $\pi_i$ by $y$. If $i\in A(\pi)$ or $i\in P(\pi)$, we label $\pi_{i-1}$ by $x$. The weight of a marked Stirling permutation $\pi$ of $[n]_2$ with $m$ descents is given by
\[
w(\pi)=x^{2n+1-m}y^{m}.
\]

\begin{theo}\label{S_grammar}
For $n\geqslant 1$, we have
 $$
 D^n(x)=\sum_{m=1}^nT(n,m)x^{2n-m+1}y^m.
 $$
 Setting $x=1$, we have
 $$
 D^n(x)|_{x=1}=T_n(y).
 $$
\end{theo}
\begin{pf}\label{}
We aim to show that ${D}^{n}(x)$ equals the sum of the weights of marked Stirling permutations of $[n]_2$.  We use induction on $n$. The case for $n=0$ follows from the fact that the weight of the empty permutation is $x$.
Assume that the theorem holds for $n-1$, that is,
\[
D^{n-1}(x)=\sum_{\pi\in \bar{Q}_{n-1}}w(\pi).
\]

 We now use an example to demonstrate the action of $D$ on a marked Stirling permutation of $[n-1]_2$.
Let $\pi=12\bar{2}331$ with the following grammatical labeling
\[
\stackrel{}{\mbox{$x$}}\ \stackrel{1}{\mbox{$x$}}\
\stackrel{2}{\mbox{$x$}}\ \stackrel{\bar{2}}{\mbox{$x$}}\
\stackrel{3}{\mbox{$x$}}\ \stackrel{3}{\mbox{$y$}}\
\stackrel{1}{\mbox{$y$}}.
\]

If we apply the substitution rule $x\rightarrow x^2y$ to
 the fourth letter $x$, then we insert the two elements $44$ after $\bar{2}$. We keep all the labels except that
 we assign the labels $x$ and $y$ to the two new letters $44$. It is not difficult to see that
the generated marked Stirling permutation  has a consistent grammatical labeling
\[
\stackrel{}{\mbox{$x$}}\ \stackrel{1}{\mbox{$x$}}\
\stackrel{2}{\mbox{$x$}}\ \stackrel{\bar{2}}{\mbox{$x$}}\
\stackrel{4}{\mbox{$x$}}\ \stackrel{4}{\mbox{$y$}}\
\stackrel{3}{\mbox{$x$}}\ \stackrel{3}{\mbox{$y$}}\
\stackrel{1}{\mbox{$y$}}.
\]

If we apply the substitution rule $y\rightarrow 2x^2y$ to the first letter $y$, then we insert $44$ after the second element $3$. We change the label of the second element $3$ from $y$ into $x$ and assign $x$ and $y$ to the two new elements $44$.  According to the marking rule, the second element $3$ may be marked or unmarked. These two choices correspond the coefficient $2$ in the substitution rule $y\rightarrow 2x^2y$. So we are led to the following
two marked Stirling permutations with consistent grammatical labelings,
\[
\stackrel{}{\mbox{$x$}}\ \stackrel{1}{\mbox{$x$}}\
\stackrel{2}{\mbox{$x$}}\ \stackrel{\bar{2}}{\mbox{$x$}}\
\stackrel{3}{\mbox{$x$}}\ \stackrel{3}{\mbox{$x$}}\
\stackrel{4}{\mbox{$x$}}\ \stackrel{4}{\mbox{$y$}}\
\stackrel{1}{\mbox{$y$}},
\]
and
\[
\stackrel{}{\mbox{$x$}}\ \stackrel{1}{\mbox{$x$}}\
\stackrel{2}{\mbox{$x$}}\ \stackrel{\bar{2}}{\mbox{$x$}}\
\stackrel{3}{\mbox{$x$}}\ \stackrel{\bar{3}}{\mbox{$x$}}\
\stackrel{4}{\mbox{$x$}}\ \stackrel{4}{\mbox{$y$}}\
\stackrel{1}{\mbox{$y$}}.
\]
It can be verified that the above process generates all marked Stirling permutations of $[n]_2$. It follows that
$$
D^n(x)=D(D^{n-1}(x))=D\left(\sum_{\pi\in \bar{Q}_{n-1}}w(\pi)\right)=\sum_{\sigma\in \bar{Q}_n}w(\sigma).
$$
Hence the proof is complete by induction.
\qed

\end{pf}

As a multivariate   refinement of $T_n(x)$,  we define the
following generating polynomial of marked Stirling permutations
of $[n]_2$,
\[
T_n(\x,\y,\z)=\sum_{\pi\in{\bar{Q}_n}}
\prod_{i\in {A}(\pi) }x_{\pi_{i}}\prod_{i\in {D}(\pi) }y_{\pi_i}\prod_{i\in {P}(\pi) }z_{\pi_i}.
\]
Let
$$
G_n=\{x_i\rightarrow x_ny_nz_n,y_i\rightarrow 2x_ny_nz_n,z_i\rightarrow x_ny_nz_n, 0\leq i<n\}.
$$
The grammatical labeling of a marked Stirling permutation
 can be described as follows. For a marked Stirling permutation $\pi$ of $[n]_2$, if $i\in A(\pi)$, we label $\pi_{i-1}$ by $x_i$; if $i\in D(\pi)$, we label $\pi_i$ by $y_i$; and if $i\in P(\pi)$, we label $\pi_{i-1}$ by $z_i$. Then the weight of $\pi$  equals
\[
w(\pi)=\prod_{i\in {A}(\pi) }x_{\pi_{i}}\prod_{i\in {D}(\pi) }y_{\pi_i}\prod_{i\in {P}(\pi) }z_{\pi_i}.
\]

The following theorem shows that the polynomials $T_n(\x,\y,\z)$ can be
generated by the grammars $G_n$.

\begin{theo} For $n\geq 1$, we have
$$D_{n}D_{n-1}\cdots D_{1}(z_{0})=T_n(\x,\y,\z).
$$
\end{theo}

The proof of the above theorem is analogous to that of Theorem \ref{S_grammar}.
Hence the details are omitted. Here we use an
example to illustrate  the action of $D_4$ to the above marked Stirling permutation $\pi=12\bar{2}331$ with the grammatical labeling
\[
\stackrel{}{\mbox{$x_1$}}\ \stackrel{1}{\mbox{$x_2$}}\
\stackrel{2}{\mbox{$z_2$}}\ \stackrel{\bar{2}}{\mbox{$x_3$}}\
\stackrel{3}{\mbox{$z_3$}}\ \stackrel{3}{\mbox{$y_3$}}\
\stackrel{1}{\mbox{$y_1$}}.
\]
If we apply the substitution rule $x_3\rightarrow x_4y_4z_4$ of $G_4$ to the letter $x_3$, then we insert the two elements $44$ after $\bar{2}$ to get a marked Stirling  permutation with the following consistent grammatical labeling
\[
\stackrel{}{\mbox{$x_1$}}\ \stackrel{1}{\mbox{$x_2$}}\
\stackrel{2}{\mbox{$z_2$}}\ \stackrel{\bar{2}}{\mbox{$x_4$}}\
\stackrel{4}{\mbox{$z_4$}}\ \stackrel{4}{\mbox{$y_4$}}\
\stackrel{3}{\mbox{$z_3$}}\ \stackrel{3}{\mbox{$y_3$}}\
\stackrel{1}{\mbox{$y_1$}}.
\]
If we apply the substitution rule $y_3\rightarrow 2x_4y_4z_4$ of $G_4$ to the letter $y_3$, then we insert $44$ after the second element $3$ to get the following two marked Stirling  permutations with consistent grammatical labelings
\[
\stackrel{}{\mbox{$x_1$}}\ \stackrel{1}{\mbox{$x_2$}}\
\stackrel{2}{\mbox{$z_2$}}\ \stackrel{\bar{2}}{\mbox{$x_3$}}\
\stackrel{3}{\mbox{$z_3$}}\ \stackrel{3}{\mbox{$x_4$}}\
\stackrel{4}{\mbox{$z_4$}}\ \stackrel{4}{\mbox{$y_4$}}\
\stackrel{1}{\mbox{$y_1$}},
\]
and
\[
\stackrel{}{\mbox{$x_1$}}\ \stackrel{1}{\mbox{$x_2$}}\
\stackrel{2}{\mbox{$z_2$}}\ \stackrel{\bar{2}}{\mbox{$x_3$}}\
\stackrel{3}{\mbox{$z_3$}}\ \stackrel{\bar{3}}{\mbox{$x_4$}}\
\stackrel{4}{\mbox{$z_4$}}\ \stackrel{4}{\mbox{$y_4$}}\
\stackrel{1}{\mbox{$y_1$}}.
\]

For $n=0$, the empty permutation is labeled by $z_0$.
For $n=1,2$,  $T_n(\x,\y,\z)$ are given below,
\begin{eqnarray*}
D_1(z_{0})&=&\stackrel{}{\mbox{$x_1$}}\
\stackrel{1}{\mbox{$z_1$}}\ \stackrel{1}{\mbox{$y_1$}},\\[5pt]
D_2D_1(z_{0})&=&
\stackrel{}{\mbox{$x_2$}}\ \stackrel{2}{\mbox{$z_2$}}\
\stackrel{2}{\mbox{$y_2$}}\ \stackrel{1}{\mbox{$z_1$}}\
\stackrel{1}{\mbox{$y_1$}}+
\stackrel{}{\mbox{$x_1$}}\ \stackrel{1}{\mbox{$x_2$}}\
\stackrel{2}{\mbox{$z_2$}}\ \stackrel{2}{\mbox{$y_2$}}\
\stackrel{1}{\mbox{$y_1$}}+\stackrel{}{\mbox{$x_1$}}\
\stackrel{1}{\mbox{$z_1$}}\ \stackrel{1}{\mbox{$x_2$}}\ \stackrel{2}{\mbox{$z_2$}}\ \stackrel{2}{\mbox{$y_2$}}\\[5pt]
&&\quad+
\stackrel{}{\mbox{$x_1$}}\ \stackrel{1}{\mbox{$z_1$}}\
\stackrel{\bar{1}}{\mbox{$x_2$}}\ \stackrel{2}{\mbox{$z_2$}}\ \stackrel{2}{\mbox{$y_2$}}.
\end{eqnarray*}

\end{section}

\begin{section}{Grammars preserving stability}\label{theory}

In this section, we
prove the stability of the multivariate polynomials $B_n(\x,\y,\z,\u,\v)$ and $T_n(\x,\y,\z)$
based on context-free grammars and the characterization of stability preserving
linear operators due to Borcea and Br\"{a}nd\'{e}n\cite{BB}.

Our idea of proving the stability of the polynomials by a sequence of context-free grammars $\{G_n\}$ goes as follows. Since the initial polynomial $x$ is stable, if $D_1,D_2,\ldots,D_n,\ldots$ preserve stability,  then  $D_nD_{n-1}\ldots D_1(x)$ is stable.
 If $D_n$ is not stability preserving, then we try to find a sequence of stability preserving
operator $\{T_n\}$  such that
 \[T_nT_{n-1}\ldots T_1(x)=D_nD_{n-1}\ldots D_1(x).
 \]
 If such operators $T_n$ exist, then we reach the
 conclusion that the multivariate
 polynomials $D_nD_{n-1}\ldots D_1(x)$ are stable.

 Note that  $B_n(\x,\y,\z,\u,\v)$ and $T_n(\x,\y,\z)$ are
all  multiaffine polynomials in the sense that
the degree in each variable is at most $1$.
In order to construct the stability preserving operators
 $T_n$ based on the grammars $G_n$, we
 consider some equivalent forms of production rules
 when we restrict our attention to multiaffine polynomials.

  For example, let
\[
G_n=\{a\rightarrow ab_n,b_i\rightarrow b_n,0\leq i<n\}.
\]
Observe that as far as the computation is concerned,
the formal differential operator $D_n$ with respect to $G_n$
is in accordance with  the following operator
\[
T_n=b_{n}(1+\sum_{i=1}^n\partial/\partial_{b_i}),
\]
 when they are applied to certain polynomials.
 Thus we obtain
\[
T_nT_{n-1}\ldots T_1(a)=D_nD_{n-1}\ldots D_1(a).
\]
However,  $D_n$ and $T_n$ are different operator in general, since
\[
D_n(a+b_1)\neq T_n(a+b_1).
\]

For multiaffine polynomials, the characterization of
stability preserving operators is simpler than  that for the
general case.
For the purpose of this paper,
we only need  the following sufficient condition
 to prove the stability of $B_n(\x,\y,\z,\u,\v)$ and $T_n(\x,\y,\z)$,
 see Borcea and Br\"{a}nd\'{e}n \cite{BB}.

\begin{lem} \label{HV}
Let $f\in \mathbb{C}[z_1,z_2,\dots,z_n]$ be a stable multiaffine polynomial and let $T$ denote a linear operator acting on the polynomials in $\mathbb{C}[z_1,z_2,\dots,z_n]$. If
\[
T\left(\prod_{i=1}^n(z_i+w_i)\right)
\in\mathbb{C}[z_1,z_2,\dots,z_n,w_1,\ldots,w_n]
\]
is stable, then $T(f)$ is either stable or identically 0.
\end{lem}

Next we show how to
prove the stability of  polynomials generated by context-free
grammars. Let us consider the multiaffine polynomials $C_n(\x,\y,\z)$ defined by Haglund and Visontai \cite{HV}.
Let
$$
G_n=\{x_i\rightarrow x_n y_n z_n,y_i\rightarrow x_n y_n z_n,z_i\rightarrow x_n y_n z_n, 0\leq i<n\},
$$
and let $G_n$ denote the differential operator associated
with the grammar $G_n$. Let $f_n=D_nD_{n-1}\cdots D_{1}(z_{0})$. From the grammatical labelings, it is clear that $f_n$ is multiaffine. We wish to prove
 the stability of $f_n$ by induction on $n$.
 Since $z_{0}$ is stable, it suffices to prove that the operator $D_{n+1}$ preserves stability of multiaffine polynomials.

Let
\[
F=\prod_{i=0}^n(x_i+w_i)(y_i+v_i)(z_i+u_i).
\]
By Lemma \ref{HV}, it  suffices to check the stability of $D_{n+1}(F)$, that is,
\[
D_{n+1}(F)=x_{n+1}y_{n+1}z_{n+1}F
\sum_{i=0}^n\left(\frac{1}{x_i+w_i}
+\frac{1}{y_i+v_i}+\frac{1}{z_i+u_i}\right)
\]
is stable.

If $x_i$, $y_i$, $z_i$, $w_i$, $v_i$ and $u_i$ have positive imaginary parts for all $0\leq i\leq n$, then
$$
\xi=\sum_{i=0}^n\left(\frac{1}{x_i+w_i}
+\frac{1}{y_i+v_i}+\frac{1}{z_i+u_i}\right)
$$
has negative imaginary part. Thus,
\[
D_{n+1}(F)=x_{n+1}y_{n+1}z_{n+1}F\xi \not= 0.
\]
Hence $D_{n+1}(F)$ is stable. By Lemma \ref{HV}, we find that $D_{n+1}(f_n)$ is a stable polynomial. So we conclude that
\[
f_{n+1}=D_{n+1}D_nD_{n-1}\cdots D_{1}(z_{0})
\]
is stable.

The stability of $A_n(\x,\y)$ can be proved in the
same way.  Indeed, let
\[
G_n=\{x_i\rightarrow x_n y_n,y_i\rightarrow x_n y_n, 0\leqslant i< n\},
\]
and let $D_n$ denote the
  differential operator  with respect to $G_n$.
  It turns out that the operator $D_n$
    preserves the stability of multiaffine polynomials.

It is worth mentioning that the formal differential operators
used in the above two examples are essentially equivalent to the operators given  by Haglund and Visontai \cite{HV} in their
proofs of the stability of $C_n(\x,\y,\z)$ and $A_n(\x,\y)$.

Next we construct stable multivariate refinements of $S_n(x)$, the Stirling polynomials. Recall that
 the grammar \[ G=\{a \rightarrow
 ab, b \rightarrow b\}
 \]
 generates the polynomials $S_n(x)$. Define
$$
G_n=\{a\rightarrow ab_n,b_i\rightarrow b_n, 1\leq i<n\},
$$
and let $D_n$ denote the formal differential operator associated
 with $G_n$.
We define the grammatical labeling of a partition as follows. For a partition $P=\{P_1,P_2,\ldots,P_k\}$, we label the partition itself by the letter $a$ and label a block $P_i$ by the letter $b_m$, where $m$ is the maximum element in $P_i$. Then the weight of $P$ is given by
\[
w(P)=a\prod_{i=1}^kb_{m_i},
\]
where $m_i$ is the maximum element in $P_i$. Denote by $S_n(a,\b)$ the sum of  weights of partitions of $[n]$. The next theorem shows that
$S_n(a,\b)$ can be generated by $G_n$. However, in this case,
the differential operator $D_n$ associated with $G_n$ is not
stability preserving even for multiaffine polynomials. Instead, we
shall find an equivalent operator $T_n$ that preserves stability
for multiaffine polynomials.

\begin{theo}
For $n\geq 1$, we have
\[
D_nD_{n-1}\ldots D_1(a)=S_n(a,\b).
\]
\end{theo}

The proof of the above theorem is analogous to that of
(\ref{partition}). Here we use the same example to
 demonstrate the action of $D_7$ on a partition of $[6]$.
Consider the  partition of $\{1, 2, 3, 4, 5, 6\}$ with the following
grammatical labeling
\[
\stackrel{\{1,3,6\}}{\mbox{$b_6$}}\stackrel{\{2,5\}}{\mbox{$b_5$}}
     \stackrel{\{4\}}{\mbox{$b_4$}}\ \stackrel{}{\mbox{$a$}}.
     \]

If we apply the substitution rule $a\rightarrow ab_7$ of $G_7$ to the letter $a$, then
 we get a partition with a consistent labeling
\[\stackrel{\{1,3,6\}}{\mbox{$b_6$}}\stackrel{\{2,5\}}{\mbox{$b_5$}} \stackrel{\{4\}}{\mbox{$b_4$}}
\stackrel{\{7\}}{\mbox{$b_7$}}\ \stackrel{}{\mbox{$a$}}.
\]

If we apply the substitution rule $b_5\rightarrow b_7$ of $G_7$ to the letter $b_5$, then
we get  the following partition with a
consistent grammatical labeling
\[\stackrel{\{1,3,6\}}{\mbox{$b_6$}}\stackrel{\{2,5,7\}}{\mbox{$b_7$}} \stackrel{\{4\}}{\mbox{$b_4$}}\
\stackrel{}{\mbox{$a$}}.
\]

\begin{theo}
For $n\geq 1$, the multivariate polynomial $S_n(a,\b)$ is stable.
\end{theo}

\begin{pf}
From the grammatical labelings, we see  that $S_n(a,\b)$ is multiaffine. Note that  $S_n(a,\b)$ is multiaffine in $a, b_1, b_2, \ldots, b_n$ with every term containing $a$ as a factor.
Since
\[ G_{n+1}=\{a\rightarrow ab_{n+1},b_i\rightarrow b_{n+1},1\leq i\leq n\},\]
for each multiaffine monomial of $S_n(a,\b)$ which is of the form
$ah$,
we have
\[ D_{n+1}(ah)=ab_{n+1} h + a D_{n+1}(h).\]
It follows that
\begin{eqnarray*}
S_{n+1}(a,\b)&=&D_{n+1}(S_n(a,\b))\\[5pt]
&=&b_{n+1}S_n(a,\b)
+b_{n+1}\sum_{i=1}^n\partial/\partial_{b_i}(S_n(a,\b)).
\end{eqnarray*}
Define
\[
T_{n+1}=b_{n+1}\left(1+\sum_{i=1}^n\partial/\partial_{b_i}\right).
\]
Then we have $S_{n+1}(a,\b)=T_{n+1}(S_n(a,\b))$.

We proceed to prove
 the stability of $S_n(a,\b)$ by induction on $n$.
Since $a$ is stable, we only need to show that the linear operator $T_{n+1}$ preserves stability of multiaffine polynomials.

Let
\[
F=(a+w)\prod_{i=1}^n(b_i+v_i).
\]
Then we have
\begin{eqnarray*}
T_{n+1}(F)&=&b_{n+1}F+b_{n+1}\sum_{i=1}^n\partial/\partial_{b_i}(F)\\[5pt]
&=&b_{n+1}F+b_{n+1}F\sum_{i=1}^n\frac{1}{b_i+v_i}\\[5pt]
&=&b_{n+1}F \left(1+\sum_{i=1}^n\frac{1}{b_i+v_i}\right),
\end{eqnarray*}
To prove that $T_{n+1}(F)$ is stable, we assume that $a$, $w$, $b_i$ and $v_i$ have positive imaginary parts for all $1\leq i\leq n+1$.
Consequently,
$$
\xi=1+\sum_{i=1}^n\frac{1}{b_i+v_i}
$$
is nonzero since it has negative imaginary part.
Moreover, each factor of $F$ has positive imaginary part, and so does $b_{n+1}$. This yields that $F$ and $b_{n+1}$ do not vanish.  It follows that
 \[ {T_{n+1}}(F) =b_{n+1}F \xi \neq 0.\]
 Hence $T_{n+1}(F)$ is stable.
In view of Lemma \ref{HV}, we see that
$T_{n+1}(S_n(a,\b))$ is stable.
This completes the proof.
\qed
\end{pf}

It is worth mentioning
that we use the operator $T_{n+1}$ instead of $D_{n+1}$
in the above proof because the operator $D_{n+1}$ does not satisfy the condition in Lemma \ref{HV}. Take $D_2$ as an example.
It can be seen that $D_2((a+w)(b_1+u))$ is not stable.
Note that
\[
D_2((a+w)(b_1+u)) = b_2(a(b_1+u+1)+w).\]
Let $a=\frac{i-1}{2},b_1=\frac{i}{2}-1,u=\frac{i}{2}-1$ and $w=i$.
But we have $D_2((a+w)(b_1+u))=0$. This implies that $D_2((a+w)(b_1+u))$ is not stable.

Next we prove the stability of $B_n(\x,\y,\z,\u,\v)$.

\begin{theo}\label{B_n}
For $n\geq 1$, the multivariate polynomial $B_n(\x, \y, \z,\u,  \v)$ is stable.
\end{theo}

\begin{pf}
 Let $f_n= D_{n}D_{n-1}\ldots D_1(x_0)$.  From the grammatical
labelings, it can be seen that $f_n$ is multiaffine.  We proceed to prove the stability of $D_{2n}D_{2n-1}\ldots D_1(x_0)$ by induction on $n$. The stability of    $x_0$ is evident.

We now assume that $f_{2n-2}$ is stable. Let us
 consider the actions of $D_{2n-1}$ and $D_{2n}$.
By using the argument in the proof of the stability of
 $C_n(\x,\y,\z)$, it can be shown that
  the operator $D_{2n-1}$ preserves stability of multiaffine polynomials. This leads to the stability of $f_{2n-1}$
  since $f_{2n-1}=D_{2n-1}(f_{2n-2})$.

Recall that
\begin{eqnarray*}
G_{2n}&=&\{x_i,y_i,z_i,u_i,v_i\rightarrow x_ny_nz_n, 0\leqslant i<n;\\[5pt]
 &&\quad \quad u_{n}\rightarrow x_{n}z_nu_{n}, v_{n}\rightarrow x_{n}y_{n}z_n \}.
\end{eqnarray*}
Let $B$ denote the following alphabet
\[
\{x_i,y_i,z_i,u_i,v_i,0\leqslant i< n\}\cup \{v_n\}.
\]
Since $f_{2n-1}$ is multiaffine and each term in $f_{2n-1}$ contains $u_n$, we may write  a monomial of
 $f_{2n-1}$ in the form $u_nh$. Then we have
\begin{eqnarray*}
D_{2n}(u_nh)&=&(x_nz_nu_n)h+x_n y_nz_n  D_{2n}(h).
\end{eqnarray*}
Thus,
\begin{eqnarray*}
f_{2n}&=&D_{2n}(f_{2n-1})\\[5pt]
&=&x_nz_nf_{2n-1}
+x_ny_nz_n\sum_{w\in B}\partial/\partial_{w}(f_{2n-1}).
\end{eqnarray*}
Hence we may write $f_{2n}$ as $T(f_{2n-1})$, where $T$ is
a linear operator as given by
\[
T=x_nz_n
+x_ny_nz_n\sum_{w\in B}\partial/\partial_{w}.
\]

It remains to show that $T$ preserves the stability of multiaffine polynomials. Let
\[
F=(u_n+r_{u_n})\prod_{w\in B}(w+r_w).
\]
By Lemma \ref{HV}, it suffices to verify  the stability of the  following polynomial
\begin{eqnarray*}
T\left(F\right)&=&x_nz_nF+x_ny_nz_nF\sum_{w\in B}\frac{1}{w+r_w}\\[5pt]
&=&x_ny_nz_nF\left(\frac{1}{y_n}+\sum_{w\in B}\frac{1}{w+r_w} \right).
\end{eqnarray*}
Suppose that all the variables $x_i,y_i,z_i,u_i,v_i,r_{x_i},r_{y_i},r_{z_i},r_{u_i}$ and $r_{v_i}$ have positive imaginary parts for $0\leq i\leq n$. Then
$$
\xi=\frac{1}{y_n}+\sum_{w\in B}\frac{1}{w+r_w}
$$
 has negative imaginary part, and so it is nonzero. Meanwhile, every factor of $F$ is nonzero since its  imaginary part is positive.
  Note that under the above assumption, $x_n$, $y_n$ and $z_n$
  have positive imaginary parts, and hence they are nonzero. Consequently, $T\left(F\right)=x_ny_nz_nF\xi$ does not vanish. This leads to the stability of  $T(F)$.

   In light of Lemma \ref{HV}, we deduce that $f_{2n}$ is stable.
   This completes the proof.
\qed
\end{pf}


The proof of the stability of $C_n(\x,\y,\z)$  applies
  to the stability of $T_n(\x,\y,\z)$. The details are omitted.

\begin{theo}\label{cgrammar}
For $n\geqslant 1$, the multivariate polynomial $T_n(\x,\y,\z)$ is stable.
\end{theo}

Multivariate stable polynomials can be reduced to real-rooted univariate   polynomials by diagonalization and specialization, see Wagner \cite{Wagner}. More precisely, if   $f\in\mathbb{R}[z_1,z_2,\ldots,z_n]$ is stable, then $f(z_1,\ldots,z_n)|_{z_i=z_j}$ and $ f(z_1,\ldots,z_n)|_{z_i=a}$ are also stable, where $1\leqslant i\neq j\leqslant n$ and $a\in\mathbb{R}$.
For example, setting $a=1$ and $b_1=b_2=\cdots=b_n=x$ in $S_n(a,\b)$ leads to the real-rootedness of $S_n(x)$, see Harper \cite{Harper}.

For the multivariate stable polynomials $B_n(\x,\y,\z,\u,\v)$, applying the diagonalization $y_i=v_i=x$ and the specialization  $x_i=u_i=z_i=1$ for all $0\leq i\leq n$, we are led to Theorem \ref{Egge}.

Let $M(n,k)$ denote the number of Legendre-Stirling permutations of $M_n$ with $k$ barred descents. By setting $x_i=u_i=y_i=z_i=1$, and $v_i=x$ for all $0\leq i\leq n$, we obtain the real-rootedness of the generating function of $M(n,k)$.

\begin{coro}
For $n\geqslant 1$, the polynomial
\[
M_n(x)=\sum_{k=1}^{n}M(n,k)x^k
\]
 has only real roots.
\end{coro}

For the multivariate stable polynomials $T_n(\x,\y,\z)$,
 by setting $x_i=z_i=1$ and $y_i=y$ for all $0\leqslant i\leqslant n$, we are led to the real-rootedness of $T_n(y)$, which is
 equivalent to the real-rootedness of $C_n(x)$.

\end{section}

\vspace{.2cm} \noindent{\bf Acknowledgments.}
We wish to thank Daniel K. Du, Cindy C.Y. Gu and Arthur L.B. Yang for helpful comments.
This work was supported
by the PCSIRT Project of the Ministry of Education, and the National
Science Foundation of China.

\end{document}